\documentclass[11pt,a4paper]{amsart}
\usepackage{amssymb}

\usepackage{amsfonts}

\newtheorem{theorem}{Theorem}
\theoremstyle{plain}
\newtheorem{acknowledgement}{Acknowledgement}

\newtheorem{corollary}{Corollary}

\newtheorem{definition}{Definition}

\newtheorem{lemma}{Lemma}

\newtheorem{proposition}{Proposition}
\newtheorem{remark}{Remark}

\numberwithin{equation}{section}
\input{tcilatex}

\begin{document}
\title[On the Digraph of a Unitary Matrix]{On the Digraph of a Unitary Matrix%
}
\author{Simone Severini}
\address{Computer Science, Univ. Bristol, Bristol, U.K.}
\email{severini@cs.bris.ac.uk}
\date{May 2002: Published in SIAM Journal on Matrix Analysis and
Applications, Volume 25, Number 1, pp. 295-300, July 2003.}
\subjclass[2000]{05C20, 51F25, 81P68 \\
DOI. 10.1137/S0895479802410293}

\begin{abstract}
Given a matrix $M$ of size $n$, the digraph $D$ on $n$ vertices is said to
be the \emph{digraph of} $M$, when $M_{ij}\neq 0$ if and only if $\left(
v_{i},v_{j}\right) $ is an arc of $D$. We give a necessary condition, called
strong quadrangularity, for a digraph to be the digraph of a unitary matrix.
With the use of such a condition, we show that a line digraph, $%
\overrightarrow{L}D$, is the digraph of a unitary matrix if and only if $D$
is Eulerian. It follows that, if $D$ is strongly connected and $%
\overrightarrow{L}D$ is the digraph of a unitary matrix then $%
\overrightarrow{L}D$ is Hamiltonian. We conclude with some elementary
observations. Among the motivations of this paper are coined quantum random
walks, and, more generally, discrete quantum evolution on digraphs. 
\end{abstract}

\maketitle

\section{Introduction}

Let $D=\left( V,A\right) $ be a digraph on $n$ vertices, with labelled
vertex set $V\left( D\right) $, arc set $A\left( D\right) $ and adjacency
matrix $M\left( D\right) $. We assume that $D$ may have loops and multiple
arcs. Let $M$ be a matrix over any field. A digraph $D$ is the \emph{digraph
of }$M$, or, equivalently, the \emph{pattern of }$M$, if $%
{\vert}%
V\left( D\right) 
{\vert}%
=n$, and, for every $v_{i},v_{j}\in V(D)$, $\left( v_{i},v_{j}\right) \in
A\left( D\right) $ if and only if $M_{ij}\neq 0$. The \emph{support} $^{s}M$
of the matrix $M$ is the $\left( 0,1\right) $-matrix with element 
\begin{equation*}
^{s}M_{ij}=\left\{ 
\begin{tabular}{ll}
1 & if $M_{ij}\neq 0,$ \\ 
0 & otherwise.%
\end{tabular}%
\ \right. 
\end{equation*}%
Then the digraph of a matrix is the digraph whose adjacency matrix is the
support of the matrix. The \emph{line digraph} of a digraph $D$, denoted by $%
\overrightarrow{L}D$, is the digraph whose vertex set $V(\overrightarrow{L}D)
$ is $A(D)$ and $(\left( v_{i},v_{j}\right) ,\left( v_{j},v_{k}\right) )\in
A(\overrightarrow{L}D)$ if and only if $(v_{i},v_{j}),(v_{j},v_{k})\in A(D)$.

A \emph{discrete quantum random walk} on a digraph $D$ is a discrete walk on 
$D$ induced by a unitary transition matrix. The term \emph{quantum random
walk} was coined by Gudder (see, \emph{e.g.}, \cite{G88}), who introduced
the model and proposed to use it to describe the motion of a quantum object
in discrete space-time and to describe the internal dynamics of elementary
particles. Recently, quantum random walks have been rediscovered, in the
context of quantum computation, by Ambainis \emph{et al.} (see \cite{ABNVW01}
and \cite{AKV01}). Since the notion of quantum random walks is analogous to
the notion of random walks, interest on quantum random walks has been
fostered by the successful use of random walks on combinatorial structures
in probabilistic algorithms (see, \emph{e.g.}, \cite{L93}). Clearly, a
quantum random walk on a digraph $D$ can be defined if and only if $D$ is
the digraph of a unitary matrix. Inspired by the work of David Meyer on
quantum cellular automata \cite{M96}, the authors of \cite{ABNVW01} and \cite%
{AKV01} overcame this obstacle in the following way.\ In order to define a
quantum random walk on a simple digraph $D$, which is regular and is not the
digraph of a unitary matrix, a quantum random walk on $\overrightarrow{L}D$
is defined. The digraph $\overrightarrow{L}D$ is the digraph of a unitary
matrix. When we chose an appropriate labeling for $V(\overrightarrow{L}D)$,
a quantum random walk on $\overrightarrow{L}D$ induces a probability
distribution on $V(D)$. The quantum random walk on $\overrightarrow{L}D$ is
called the \emph{coined quantum random walk} on $D$.

With this scenario in mind, the question which this paper addresses is the
following: On which digraphs can quantum random walks be defined? In a more
general language, we are interested in the combinatorial properties of the
digraphs of unitary matrices. We give a simple necessary condition, called 
\emph{strong quadrangularity}, for a digraph to be digraph of a unitary
matrix. While it seems too daring to conjecture that such a condition is
sufficient in the general case, we discover \textquotedblleft
accidentally\textquotedblright\ that strong quadrangularity is sufficient
when the digraph is a line digraph. We also prove that if a line digraph of
a strongly connected digraph is the digraph of a unitary matrix, then it is
Hamiltonian. We observe that strong quadrangularity is sufficient to show
that certain strongly regular graphs are digraphs of unitary matrices and
that $n$-paths, $n$-paths with loops at each vertex, $n$-cycles, directed
trees and trees are not. In \cite{GZe88} and \cite{M96} the fact that an $n$%
-path is not the digraph of a unitary matrix was called the \emph{NO-GO Lemma%
}. A consequence of the lemma was that there is no nontrivial, homogeneous,
local, one-dimensional quantum cellular automaton. Proposition \ref{mey}
below can be then interpreted as a simple combinatorial version of the NO-GO
Lemma.

We refer to \cite{T84} and to \cite{BR91}, for notions of graph theory and
matrix theory, respectively.

\section{Digraphs of unitary matrices}

Let $D=(V,A)$ be a digraph. A vertex of a digraph is called \emph{source} (%
\emph{sink}) if it has no ingoing (outgoing) arcs. A vertex of a digraph is
said to be \emph{isolated} if it is not joined to another vertex. We assume
that $D$ has no sources, sinks and disconnected loopless vertices. By this
assumption, $A(D)$ has neither zero-rows nor zero-columns. For every $%
S\subset V(D)$, denote by 
\begin{equation*}
\begin{tabular}{ccc}
$N^{+}\left[ S\right] =\{v_{j}:(v_{i},v_{j})\in A(D),v_{i}\in S\}$ & and & $%
N^{-}\left[ S\right] =\left\{ v_{i}:(v_{i},v_{j})\in A(D),v_{j}\in S\right\} 
$%
\end{tabular}%
\end{equation*}
the \emph{out-neighbourhood} and \emph{in-neighbourhood} of $S$,
respectively. Denote by $%
{\vert}%
X%
{\vert}%
$ the cardinality of a set $X$. The non-negative integers $%
{\vert}%
N^{-}\left[ v_{i}\right] 
{\vert}%
$ and $%
{\vert}%
N^{+}\left[ v_{i}\right] 
{\vert}%
$ are called \emph{invalency} and \emph{outvalency} of the vertex $v_{i}$,
respectively. A digraph $D$ is \emph{Eulerian} if and only if every vertex
of $D$ has equal invalency and outvalency.

The notion defined in Definition 1 is standard in combinatorial matrix
theory (see, \emph{e.g.}, \cite{BR91}). In graph theory, the term \emph{%
quadrangular} was first used in \cite{GZ98}.

\begin{definition}
A digraph $D$ is said to be \emph{quadrangular} if, for any two distinct
vertices $v_{i},v_{j}\in V(D)$, we have 
\begin{equation*}
\begin{tabular}{lll}
$\left\vert N^{+}\left[ v_{i}\right] \cap N^{+}\left[ v_{j}\right]
\right\vert \neq 1$ & and & $\left\vert N^{-}\left[ v_{i}\right] \cap N^{-}%
\left[ v_{j}\right] \right\vert \neq 1$.%
\end{tabular}%
\ 
\end{equation*}
\end{definition}

\begin{definition}
\label{squad}A digraph $D$ is said to be \emph{strongly quadrangular} if
there does not exist a set $S\subseteq V\left( D\right) $ such that, for any
two distinct vertices $v_{i},v_{j}\in S$, 
\begin{equation*}
\begin{tabular}{lll}
$N^{+}\left[ v_{i}\right] \cap \bigcup_{j\neq i}N^{+}\left[ v_{j}\right]
\neq \emptyset $ & and & $N^{+}\left[ v_{i}\right] \cap N^{+}\left[ v_{j}%
\right] \subseteq T,$%
\end{tabular}%
\ 
\end{equation*}%
where $\left\vert T\right\vert <\left\vert S\right\vert $, and similarly for
the in-neighbourhoods.
\end{definition}

\begin{remark}
\emph{Note that if a digraph is strongly quadrangular then it is
quadrangular.}
\end{remark}

\begin{lemma}
\label{sq}Let $D$ be a digraph. If $D$ is the digraph of a unitary matrix
then $D$ is strongly quadrangular.
\end{lemma}

\begin{proof}
Suppose that $D$ is the digraph of a unitary matrix $U$ and that $D$ is not
strongly quadrangular. Then there is a set $S\subseteq V\left( D\right) $
such that, for any two distinct vertices $v_{i},v_{j}\in S$, $N^{+}\left[
v_{i}\right] \cap \bigcup_{j\neq i}N^{+}\left[ v_{j}\right] \neq \emptyset $
and $N^{+}\left[ v_{i}\right] \cap N^{+}\left[ v_{j}\right] \subseteq T$
where $\left\vert T\right\vert <\left\vert S\right\vert $. This implies that
in $U$, there is a set $S^{\prime }$ of rows which contribute, with at least
one nonzero entry, to the inner product with some other rows in $S^{\prime }$%
. In addition, the nonzero entries of any two distinct rows in $S^{\prime }$%
, which contribute to the inner product of the two rows, are in the columns
of the same set of columns $T^{\prime }$ such that $\left\vert T^{\prime
}\right\vert <\left\vert S^{\prime }\right\vert $. Then the rows of $%
S^{\prime }$ form a set of orthonormal vectors of dimension smaller than the
cardinality of the set itself. This contradicts the hypothesis. The same
reasoning holds for the columns of $U$.
\end{proof}

Two digraphs $D$ and $D^{\prime }$ are \emph{permutation equivalent} if
there are permutation matrices $P$ and $Q$, such that $M\left( D^{\prime
}\right) =PM\left( D\right) Q$ (and hence also $P^{-1}M(D^{\prime
})Q^{-1}=M(D)$). If $Q=P^{-1}$, then $D$ and $D^{\prime }$ are said to be 
\emph{isomorphic}. We write $D\cong D^{\prime }$ if $D$ and $D^{\prime }$
are isomorphic. Denote by $I_{n}$ the identity matrix of size $n$. Denote by 
$A^{\intercal }$ the transpose of a matrix $A$.

\begin{lemma}
\label{equi}Let $D$ and $D^{\prime }$ be permutation equivalent digraphs.
Then $D$ is the digraph of a unitary matrix if and only if $D^{\prime }$ is.
\end{lemma}

\begin{proof}
Suppose that $D$ is the digraph of a unitary matrix $U$. Then, for
permutation matrices $P$ and $Q$, we have $PUQ=U^{\prime }$, where $%
U^{\prime }$ is a unitary matrix of the digraph $D^{\prime }$. The converse
is similar.
\end{proof}

\begin{lemma}
\label{pieni}For any $n$ the complete digraph is the digraph of a unitary
matrix.
\end{lemma}

\begin{proof}
The lemma just means that for every $n$ there is a unitary matrix without
zero entries. An example is given by the Fourier transform on the group $%
\mathbb{Z}/n\mathbb{Z}$ (see, \emph{e.g.} \cite{T99}).
\end{proof}

A digraph $D$ is said to be $\left( k,l\right) $\emph{-regular} if, for
every $v_{i}\in V\left( D\right) $, $\left| N^{-}\left[ v_{i}\right] \right|
=k$ and $\left| N^{+}\left[ v_{i}\right] \right| =l$. If $k=l$ then $D$ is
said to be simply $k$\emph{-regular}.

\begin{remark}
\label{triangolo}\emph{Not every }$k$\emph{-regular digraph is the digraph
of a unitary matrix. Let} 
\begin{equation*}
M\left( D\right) =\left[ 
\begin{array}{ccc}
0 & 1 & 1 \\ 
1 & 0 & 1 \\ 
1 & 1 & 0%
\end{array}
\right] . 
\end{equation*}
\emph{Note that }$D$\emph{\ is }$2$\emph{-regular and it is not quadrangular.%
}
\end{remark}

\begin{remark}
\emph{Not every quadrangular digraph is the digraph of a unitary matrix. Let 
} 
\begin{equation*}
M\left( D\right) =\left[ 
\begin{array}{cccc}
1 & 1 & 1 & 1 \\ 
1 & 1 & 1 & 1 \\ 
1 & 1 & 0 & 0 \\ 
1 & 1 & 0 & 0%
\end{array}
\right] . 
\end{equation*}
\emph{Note that }$D$\emph{\ is quadrangular and is not the digraph of a
unitary matrix. In fact, }$D$\emph{\ is not strongly quadrangular.}
\end{remark}

\begin{definition}
A digraph $D$ is said to be \emph{specular} when, for any two distinct
vertices $v_{i},v_{j}\in V(D)$, if $N^{+}\left[ v_{i}\right] \cap N^{+}\left[
v_{j}\right] \neq \emptyset $, then $N^{+}\left[ v_{i}\right] =N^{+}\left[
v_{j}\right] $, and, equivalently, if $N^{-}\left[ v_{i}\right] \cap N^{-}%
\left[ v_{j}\right] \neq \emptyset $ then $N^{-}\left[ v_{i}\right] =N^{-}%
\left[ v_{j}\right] $.
\end{definition}

\begin{definition}
A $n\times m$ matrix $M$ is said to have \emph{independent submatrices} $%
M_{1}$ and $M_{2}$ when, for every $1\leq i,k\leq n$ and $1\leq j,l\leq m$,
if $M_{ij}\neq 0$ is an entry of $M_{1}$ and $M_{kl}\neq 0$ is an entry of $%
M_{2}$ then $i\neq k$ and $j\neq l$.
\end{definition}

\begin{theorem}
\label{ladder}A specular and strongly quadrangular digraph is the digraph of
a unitary matrix.
\end{theorem}

\begin{proof}
Let $D$ be a digraph. Note that if $D$ is specular and strongly quadrangular
then $M\left( D\right) $ is composed of independent matrices. The theorem
follows then from Lemma \ref{pieni}.
\end{proof}

The following theorem collects some classic results on line digraphs (see, 
\emph{e.g.}, \cite{P96}).

\begin{theorem}
\label{ri}Let $D$ be a digraph.

\begin{itemize}
\item[(i)] Then, for every $\left( v_{i},v_{j}\right) \in V\left( 
\overrightarrow{L}D\right) $, 
\begin{equation*}
N^{+}\left[ \left( v_{i},v_{j}\right) \right] =N^{+}\left[ v_{j}\right] 
\text{ and }N^{-}\left[ \left( v_{i},v_{j}\right) \right] =N^{-}\left[ v_{i}%
\right] . 
\end{equation*}

\item[(ii)] A digraph $D$ is a line digraph if and only $D$ is specular.

\item[(iii)] Let $D$ be a strongly connected digraph. Then $D$ is Eulerian
if and only if $\overrightarrow{L}D$ is Hamiltonian.
\end{itemize}
\end{theorem}

\begin{corollary}
\label{main1}A strongly quadrangular line digraph is the digraph of a
unitary matrix.
\end{corollary}

\begin{proof}
The proof is obtained by point (i) of Theorem \ref{ri} together with Theorem %
\ref{ladder}.
\end{proof}

\begin{remark}
\emph{Not every line digraph which is the digraph of a unitary matrix is
Eulerian. Let } 
\begin{equation*}
\begin{tabular}{lll}
$M\left( D\right) =\left[ 
\begin{array}{cc}
1 & 1 \\ 
1 & 0%
\end{array}%
\right] $ & and & $M(\overrightarrow{L}D)=\left[ 
\begin{array}{ccc}
0 & 0 & 1 \\ 
1 & 1 & 0 \\ 
1 & 1 & 0%
\end{array}%
\right] .$%
\end{tabular}%
\ 
\end{equation*}%
\emph{Note that }$\overrightarrow{L}D$\emph{\ is not Eulerian.}
\end{remark}

In a digraph, a \emph{directed path of length }$r$, from $v_{1}$ to $v_{r+1}$%
, is a sequence of arcs of the form $\left( v_{1},v_{2}\right) ,\left(
v_{2},v_{3}\right) ,...,\left( v_{r},v_{r+1}\right) $, where all vertices
are distinct. A directed path is an \emph{Hamiltonian path} if it included
all vertices of the digraph. A directed path, in which $v_{1}=v_{r+1}$, is
called \emph{directed cycle}. An Hamiltonian path, in which $%
v_{1}=v_{r+1}=v_{n}$ and $\left\vert V\left( D\right) \right\vert =n$, is
called \emph{Hamiltonian cycle}. A digraph with an Hamiltonian cycle is said
to be \emph{Hamiltonian}.

\begin{theorem}
\label{as}Let $D$ be a digraph. Then $\overrightarrow{L}D$ is the digraph of
a unitary matrix if and only if $D$ is Eulerian or the disjoint union of
Eulerian components.
\end{theorem}

\begin{proof}
Suppose that $\overrightarrow{L}D$ is the digraph of a unitary matrix. By
Corollary \ref{main1}, $\overrightarrow{L}D$ is strongly quadrangular. If
there is $v_{i}\in V(\overrightarrow{L}D)$ such that $\left\vert N^{+}\left[
v_{i}\right] \right\vert =1$ then for every $v_{j}\in V(\overrightarrow{L}D)$%
, $N^{+}\left[ v_{i}\right] \cap N^{+}\left[ v_{j}\right] =\emptyset $.
Suppose that, for every $v_{i}\in V(\overrightarrow{L}D)$, $\left\vert N^{+}%
\left[ v_{i}\right] \right\vert =1$. Since $\overrightarrow{L}D$ is strongly
quadrangular then $A\left( D\right) =A(\overrightarrow{L}D)$ and it is a
permutation matrix. In general, for every $v_{i}\in V(\overrightarrow{L}D)$,
if $\left\vert N^{+}\left[ v_{i}\right] \right\vert =k>1$, then there is a
set $S\subset V(\overrightarrow{L}D)$ with $\left\vert S\right\vert =k-1$
and not including $v_{i}$ such that, for every $v_{j}\in S$, $N^{+}\left[
v_{j}\right] =N^{+}\left[ v_{i}\right] $. Writing $v_{i}=uv$, where $u,v\in
V\left( D\right) $, by Theorem \ref{ri}, $N^{+}\left[ v_{i}\right] =N^{+}%
\left[ v\right] $. It follows that $\left\vert N^{+}\left[ v\right]
\right\vert =k$. Then, because of $S$, it is easy to see that in $A\left(
D\right) $ there are $k$ arcs with head $w$. Hence $\left\vert N^{+}\left[ v%
\right] \right\vert =\left\vert N^{-}\left[ v\right] \right\vert $, and $D$
is Eulerian. The proof of the sufficiency is immediate.
\end{proof}

\begin{corollary}
Let $D$ be a strongly connected digraph. Let $\overrightarrow{L}D$ be the
digraph of a unitary matrix. Then $\overrightarrow{L}D$ is Hamiltonian.
\end{corollary}

\begin{proof}
We obtain the proof by point (iii)\ of Theorem \ref{ri} together with
Theorem \ref{as}.
\end{proof}

Let $G$ be a group with generating set $S$. The \emph{Cayley digraph} of $G$
in respect to $S$ is the digraph denoted by $Cay\left( G,S\right) $, with
vertex set $G$ and arc set including $\left( g,h\right) $ if and only if
there is a generator $s\in S$ such that $gs=h$.

\begin{corollary}
The line digraph of a Cayley digraph is the digraph of a unitary matrix.
\end{corollary}

\begin{proof}
The corollary follows from Theorem \ref{as}, since a Cayley digraph is
regular.
\end{proof}

A \emph{strongly regular graph} on $n$ vertices is denoted by $srg\left(
n,k,\lambda ,\mu \right) $ and is a $k$-regular graph on $n$ vertices, in
which (1) two vertices are adjacent if and only if they have exactly $%
\lambda $ common neighbours and (2) two vertices are nonadjacent if and only
if they have exactly $\mu $ common neighbours (see, \emph{e.g.}, \cite{CvL91}%
). The parameters of $srg\left( n,k,\lambda ,\mu \right) $ satisfy the
following equation: $k\left( k-\lambda -1\right) =\left( n-k-1\right) \mu $.
The disjoint union of $r$ complete graphs each on $m$ vertices, with $r,m>1$%
, is denoted by $rK_{m}$. If $m=2$ then $rK_{2}$ is called \emph{ladder graph%
}. A strongly regular graph is disconnected if and only if it is isomorphic
to $rK_{m}$.

\begin{remark}
\emph{Not every strongly regular graphs is the digraph of a unitary matrix.
The graph }$srg\left( 10,3,0,1\right) $\emph{\ is called }Petersen's graph%
\emph{. It is easy to check that }$srg\left( 10,3,0,1\right) $\emph{\ is not
quadrangular.}
\end{remark}

\begin{remark}
\emph{By Theorem \ref{ladder}, if a digraph }$D$\emph{\ is permutation
equivalent to a disconnected strongly quadrangular graph, then }$D$\emph{\
is the digraph of a unitary matrix.}
\end{remark}

The \emph{complement} of a digraph $D$ is a digraph denoted by $\overline{D}$
with the same vertex set of $D$ and with two vertices adjacent if and only
if the vertices not adjacent in $D$. A digraph $D$ is \emph{%
self-complementary} if $D\cong \overline{D}$.

\begin{remark}
\emph{The fact that }$D$\emph{\ is the digraph of a unitary matrix does not
imply that }$\overline{D}$\emph{\ is. The digraph used in the proof of
Proposition \ref{triangolo} provides a counterexample. Note that this does
not hold in the case where }$D$\emph{\ is self-complementary.}
\end{remark}

A digraph $D$ is an $n$-\emph{path}, if $V\left( D\right) =\left\{
v_{1},v_{2},...,v_{n}\right\} $ and 
\begin{equation*}
A\left( D\right) =\left\{ \left( v_{1},v_{2}\right) ,\left(
v_{2},v_{1}\right) ,\left( v_{2},v_{3}\right) ,\left( v_{3},v_{2}\right)
,...,\left( v_{n-1},v_{n}\right) ,\left( v_{n},v_{n-1}\right) \right\} ,
\end{equation*}%
where all the vertices are distinct. An $n$-path, in which $v_{1}=v_{n}$, is
called $n$\emph{-cycle}. A digraph $D$ is a \emph{directed} $n$\emph{-cycle}
if $A\left( D\right) =\left\{ \left( v_{1},v_{2}\right) ,\left(
v_{2},v_{3}\right) ,...,\left( v_{n-1},v_{1}\right) \right\} $. A digraph
without directed cycles if a \emph{directed tree}. A graph without cycle is
a \emph{tree}. 

\begin{proposition}
\label{mey}Let $D$ be a digraph. If $D$ is permutation equivalent to an $n$%
-path then it is not the digraph of a unitary matrix.
\end{proposition}

\begin{proof}
A digraph is strongly connected if and only if it is the digraph of an
irreducible matrix. Since an $n$-path is strongly connected, it is the
digraph of an irreducible matrix. Note that the number of arcs of an $n$%
-path is $2\left( n-1\right) $. The proposition is proved by Lemma \ref{sq},
together with the following result (see, \emph{e.g.}, \cite{BR91}). Let $M$
be an irreducible matrix of size $n$ and with exactly $2\left( n-1\right) $
nonzero entries. Then there is a permutation matrix $P$, such that 
\begin{equation*}
PMP^{\intercal }=\left[ 
\begin{array}{ccccc}
a_{11} & 0 & \cdots & 0 & 1 \\ 
1 & a_{22} & \cdots & 0 & 0 \\ 
\vdots & 1 & \ddots & \vdots & \vdots \\ 
0 & 0 & \cdots & \ddots & 0 \\ 
0 & 0 & \cdots & 1 & a_{nn}%
\end{array}
\right] , 
\end{equation*}
where $a_{ii}$ can be equal to zero or one. It is easy to see that for any
choice of the diagonal entries the digraph of $PMP^{\intercal }$ is not
quadrangular.
\end{proof}

\begin{proposition}
If a digraph $D$ is permutation equivalent to one of the following digraphs,
then $D$ is not the digraph of a unitary matrix: $n$-path with a loop at
each vertex, $n$-cycle, directed tree, tree.
\end{proposition}

\begin{proof}
Chosen any labeling of $D$, the proposition follows from Lemma \ref{sq} and
Lemma \ref{equi}.
\end{proof}

\begin{acknowledgement}
The author thanks Peter Cameron, Richard Jozsa, Gregor Tanner and Andreas
Winter for their help. The author is supported by a University of Bristol
research scholarship. 
\end{acknowledgement}

\end{document}